\documentclass[12pt]{article}
\usepackage{graphicx}
\usepackage{amsmath,amsthm,amssymb,enumerate}
\usepackage{euscript,mathrsfs}
\usepackage{color}
\usepackage{dsfont}
\usepackage[left=2cm,right=2cm,top=3.5cm,bottom=3.5cm]{geometry}
\usepackage{color}
\usepackage[framemethod=tikz]{mdframed}
\allowdisplaybreaks

\usepackage{soul}

\catcode`\@=11 \@addtoreset{equation}{section}

\catcode`\@=12

\newtheorem{Theorem}{Theorem}[section]
\newtheorem{Proposition}[Theorem]{Proposition}
\newtheorem{Lemma}[Theorem]{Lemma}
\newtheorem{Corollary}[Theorem]{Corollary}

\theoremstyle{definition}
\newtheorem{Definition}[Theorem]{Definition}

\newtheorem{Remark}[Theorem]{Remark}

\newcommand{\bTheorem}[1]{
	\begin{Theorem} \label{T#1} }
	\newcommand{\eT}{\end{Theorem}}

\newcommand{\bProposition}[1]{
	\begin{Proposition} \label{P#1}}
	\newcommand{\eP}{\end{Proposition}}

\newcommand{\bLemma}[1]{
	\begin{Lemma} \label{L#1} }
	\newcommand{\eL}{\end{Lemma}}

\newcommand{\bCorollary}[1]{
	\begin{Corollary} \label{C#1} }
	\newcommand{\eC}{\end{Corollary}}

\newcommand{\bRemark}[1]{
	\begin{Remark} \label{R#1} }
	\newcommand{\eR}{\end{Remark}}

\newcommand{\bDefinition}[1]{
	\begin{Definition} \label{D#1} }
	\newcommand{\eD}{\end{Definition}}

\newcommand{\wtZ}{\widetilde{Z}}

\newcommand{\Del}{\Delta_x}

\newcommand{\tvs}{\widetilde{s}}

\newcommand{\vy}{\vc{y}}

\newcommand{\vrB}{\vr_B}

\newcommand{\Ds}{\mathbb{D}_x}

\newcommand{\vuB}{\vc{u}_B}

\newcommand{\bfphi}{\boldsymbol{\varphi}}

\newcommand{\bFormula}[1]{
	\begin{equation} \label{#1}}
	\newcommand{\eF}{\end{equation}}

\newcommand{\Ov}[1]{\overline{#1}}

\newcommand{\aleq}{\stackrel{<}{\sim}}

\newcommand{\ageq}{\stackrel{>}{\sim}}

\newcommand{\vr}{\varrho}
\newcommand{\vre}{\vr_\ep}

\newcommand{\vte}{\vt_\ep}
\newcommand{\vue}{\vu_\ep}
\newcommand{\tvr}{\widetilde{\vr}}
\newcommand{\tvu}{{\tilde \vu}}
\newcommand{\tvt}{\tilde \vt}

\newcommand{\vt}{\vartheta}
\newcommand{\vu}{\vc{u}}
\newcommand{\vm}{\vc{m}}

\newcommand{\vc}[1]{{\bf #1}}

\newcommand{\Div}{{\rm div}_x}
\newcommand{\Grad}{\nabla_x}

\newcommand{\dx}{\,{\rm d} {x}}

\newcommand{\dt}{\,{\rm d} t }

\newcommand{\intO}[1]{\int_{\Omega} #1 \ \dx}

\newcommand{\D}{{\rm d}}

\newcommand{\ep}{\varepsilon}

\newcommand{\vtB}{\vt_B}

\newcommand{\br}{ \nonumber \\ }

\def\softd{{\leavevmode\setbox1=\hbox{d}%
		\hbox to 1.05\wd1{d\kern-0.4ex{\char039}\hss}}}
\definecolor{Cgrey}{rgb}{0.85,0.85,0.85}
\definecolor{Cblue}{rgb}{0.50,0.85,0.85}
\definecolor{Cred}{rgb}{1,0,0}
\definecolor{fancy}{rgb}{0.10,0.85,0.10}

\newcommand\Cbox[2]{%
	\newbox\contentbox%
	\newbox\bkgdbox%
	\setbox\contentbox\hbox to \hsize{%
		\vtop{
			\kern\columnsep
			\hbox to \hsize{%
				\kern\columnsep%
				\advance\hsize by -2\columnsep%
				\setlength{\textwidth}{\hsize}%
				\vbox{
					\parskip=\baselineskip
					\parindent=0bp
					#2
				}%
				\kern\columnsep%
			}%
			\kern\columnsep%
		}%
	}%
	\setbox\bkgdbox\vbox{
		\color{#1}
		\hrule width  \wd\contentbox %
		height \ht\contentbox %
		depth  \dp\contentbox
		\color{black}
	}%
	\wd\bkgdbox=0bp%
	\vbox{\hbox to \hsize{\box\bkgdbox\box\contentbox}}%
	\vskip\baselineskip%
}

\mdfdefinestyle{MyFrame}{%
	linecolor=black,
	outerlinewidth=1pt,
	roundcorner=5pt,
	innertopmargin=\baselineskip,
	innerbottommargin=\baselineskip,
	innerrightmargin=10pt,
	innerleftmargin=10pt,
	backgroundcolor=white!20!white}

\begin{document}


\title{Stability of planar rarefaction waves
in the vanishing dissipation limit of the Navier--Stokes--Fourier system}

\author{Eduard Feireisl	
\thanks{The work of E.F. was supported by the
		Czech Sciences Foundation (GA\v CR), Grant Agreement 24-11034S.
		The Institute of Mathematics of the Academy of Sciences of
		the Czech Republic is supported by RVO:67985840.} 
and Wladimir Neves
\thanks{W. Neves has received research grants from CNPq
	through the grants  308064/2019-4, 406460/2023-0, and also by FAPERJ 
	(Cientista do Nosso Estado) through the grant E-26/201.139/2021.
	}
}

\date{}

\maketitle

\medskip

\centerline{Institute of Mathematics of the Academy of Sciences of the Czech Republic}

\centerline{\v Zitn\' a 25, CZ-115 67 Praha 1, Czech Republic}

\bigskip

\centerline{Instituto de Matem\' atica, Universidade Federal do Rio de Janeiro}

\centerline{C.P. 68530, Cidade Universit\' aria 21945-970, Rio de Janeiro, Brazil}

\begin{abstract}
	
We consider the vanishing dissipation limit of the compressible Navier--Stokes--Fourier system, where the initial data approach 
a profile generating a planar rarefaction wave for the limit 
Euler system. We show that the associated weak solutions converge unconditionally to the planar rarefaction wave strongly in the energy norm.

\end{abstract}


{\bf Keywords:} Navier--Stokes--Fourier system, vanishing diffusion limit, planar rarefaction wave, stability. 


\section{Introduction}
\label{E}

Recently, Li, Wang, and Wang \cite{LiWanWan} showed local stability of planar rarefaction waves in the vanishing dissipation limit of the Navier--Stokes--Fourier system. 
The goal of the present paper is to extend this result to the class of global--in--time weak solutions of the same problem in the same low dissipation regime. The convergence is unconditional 
in the sense that the distance to the limit profile is estimated only in terms of the Bregman distance of the initial data evaluated by the associated relative energy.

\subsection{Navier--Stokes--Fourier system}

The time evolution of the mass density $\vr = \vr(t,x)$, the absolute temperature $\vt = \vt(t,x)$, and the velocity $\vu = \vu(t,x)$ of a general viscous, compressible, and heat conducting fluid
is governed by the \emph{Navier--Stokes--Fourier (NSF) system} of partial differential equations:
\begin{align} 
	\partial_t \vr + \Div (\vr \vu) &= 0 , \label{E1} \\
	\partial_t (\vr \vu) + \Div (\vr \vu \otimes \vu) + \Grad p &= \Div \mathbb{S},\ \label{E2}\\
	\partial_t (\vr e) + \Div (\vr e \vu) + \Grad \vc{q} &= \mathbb{S}: \Ds \vu - p \Div \vu,\ \Ds \vu \equiv \frac{1}{2} \left( \Grad \vu + \Grad^t \vu \right). \label{E3}
	\end{align} 
We consider Newtonian fluid in the zero dissipation regime. Accordingly, the viscous stress tensor $\mathbb{S}$ is given by the scaled Newton rheological law 
\begin{equation} \label{E4}
	\mathbb{S} = \mathbb{S}_\ep (\vt, \Ds \vu) = \ep \mu(\vt) \left( \Grad \vu + \Grad \vu^t - \frac{2}{d} \Div \vu \mathbb{I} \right) + 
	\ep \eta (\vt) \Div \vu \mathbb{I},\ d = 2,3,\ \ep > 0. 
	\end{equation}
Similarly, the heat flux $\vc{q}$ is determined by the Fourier law 
\begin{equation} \label{E5}
	\vc{q} = \vc{q}_\ep (\vt, \Grad \vt) = - \ep \kappa (\vt) \Grad \vt. 
\end{equation}

\subsection{Planar rarefaction wave solutions}
\label{RWS}

Motivated by the theory developed in \cite{FeiNovOpen}, we consider the mono--atomic gas equation of state (EOS)
\begin{equation} \label{E6}
	p = \frac{2}{3} \vr e.
	\end{equation}
Relation \eqref{E6} is compatible with the Second law of thermodynamics, specifically yields the existence of (specific) entropy if 
\begin{equation} \label{E7}
	p(\vr, \vt) = \vt^{\frac{5}{2}} P(Z),\ Z = \frac{\vr}{\vt^{\frac{3}{2}}}.
\end{equation} 
The associated entropy $s$ then reads 
\begin{equation} \label{E8} 
	s(\vr, \vt) = S(Z),\ S'(Z)  = -\frac{3}{2} \frac{ \frac{5}{3} P(Z) - P'(Z) Z}{Z^2},\ Z = \frac{\vr}{\vt^{\frac{3}{2}}}, 
	\end{equation}
cf. \cite[Chapter 2]{FeNo6A}.
The best known example of EOS in the form \eqref{E7} is the standard Boyle--Mariotte law, 
\begin{equation} \label{BM}
	P(Z) = Z,\ p = \vr \vt,\ e = \frac{3}{2} \vt,\ s = \frac{3}{2} \log (\vt) - \log(\vr).
\end{equation}

Consider the 1-d Riemann problem for the Euler system 
\begin{align}
\partial_t \vr + \partial_{x_1} (\vr u) &= 0, \br 
\partial_t (\vr u) + \partial_{x_1} (\vr u^2) + \partial_{x_1} p &= 0, \br
\partial_t \left( \frac{1}{2} \vr |\vu|^2 + \vr e \right) + \partial_{x_1}\left[ \left( \frac{1}{2} \vr |\vu|^2 + \vr e + p\right)  u \right] &= 0
\label{E9}
\end{align}
for $t > 0$, $x_1 \in R$, endowed with the piecewise constant data
\begin{align}
\vr(0,x_1) &= \left\{  \begin{array}{l} \tvr_L \ \mbox{for}\ x_1 < 0, \\  \tvr_R \ \mbox{for}\ x_1 \geq 0 \end{array} \right.,\ 
\vt(0,x_1) = \left\{  \begin{array}{l} \tvt_L \ \mbox{for}\ x_1 < 0, \\  \tvt_R \ \mbox{for}\ x_1 \geq 0 \end{array} \right. \br 
u(0,x_1) &= \left\{  \begin{array}{l} \widetilde{u}_L \ \mbox{for}\ x_1 < 0, \\  \widetilde{u}_R \ \mbox{for}\ x_1 \geq 0. \end{array} \right. 
\label{E10}
\end{align}

It is well known, see e.g. Chang and Hsiao \cite{CHHS}, that the Riemann problem \eqref{E9}, \eqref{E10} admits a weak solution depending only on the self--similar variable 
$\xi = \frac{x_1}{t}$. Here, we focus on the \emph{rarefaction wave} solution of \eqref{E9}, \eqref{E10} that are \emph{Lipshitz continuous} on any interal $[\delta, T]$, 
$0 < \delta < T$. The entropy $s$ associated to a rarefaction solution is necessarily constant. If, in addition, the thermodynamic function $p$, $e$ satisfy the Boyle Mariotte law 
\eqref{BM}, the rarefaction wave solutions $(\tvr, \tvt, \widetilde{u})$ coincide with their counterparts solving the isentropic Euler system, 
\begin{align} 
\tvt &= \tvr^{\frac{2}{3}} \exp \left( \frac{2}{3} s \right), \br
\partial_t \tvr &+ \partial_{x_1} (\tvr \widetilde{u}) = 0, \br
\partial_t (\tvr \widetilde{u}) &+ \partial_{x_1} (\tvr \widetilde{u}^2) + \exp \left( \frac{2}{3} s \right) \partial_{x_1} \tvr^{\frac{5}{3}} = 0.	
	\label{E11}
	\end{align}

As shown by Chen and Chen \cite{CheChe}, the planar rarefaction waves 
are unique in the class of all admissible weak solutions emanating from the same initial data even if embedded into a higher dimensional physical space.

\subsection{Stability of planar rarefaction waves}

We consider the spatial domain 
\begin{equation} \label{b1}
	\Omega = [-L,L] \times  \mathbb{T}^{d-1},  \ d = 2,3,	
	\end{equation}	
where $\mathbb{T}^{d-1}$ is the flat torus, meaning all functions defined on $\Omega$ are periodic with respect to the variables $\vy = (x_{d-1},x_d)$.  
The rarefaction waves solutions $(\tvr, \tvt, \tvu)$ introduced in the preceding section are extended to be constant with respect to $\vy$, formally,
\[
(\tvr, \tvt) (x_1, \vy) = (\tvr, \tvt) (x_1), \ 
\tvu(x_1, \vy)  = (\widetilde{u}, 0) (x_1)\ \mbox{for} \ x_1 \in (-L,L), \ \vy \in \mathbb{T}^{d-1}.
\]
As the Euler system admits a finite spead of propagation, we can choose $L = L(T) > 0$ large enough so the $(\tvr, \tvt, \tvu)$ is a solution of the Euler system 
\begin{align} 
	\partial_t \tvr + \Div (\tvr \tvu) &= 0 , \label{Eu1} \\
	\partial_t (\tvr \tvu) + \Div (\tvr \tvu \otimes \tvu) + \Grad p (\tvr, \tvt) &= 0,\ \label{Eu2}\\
	\partial_t (\tvr e (\tvr, \tvt)) + \Div (\tvr e (\tvr, \tvt) \tvu) &=  - p(\tvr, \tvt) \Div \tvu\label{Eu3}
\end{align} 
in $(0,T) \times \Omega$ satisfying the boundary conditions
\begin{align}
\tvr(t, -L, \vy) &=  \tvr_L
,\ \tvr(t, L, \vy) =  \tvr_R, \br 
\tvt(t, -L, \vy) &=  \tvt_L
,\ \tvt(t, L, \vy) =  \tvt_R \br 
\tvu(t, - L, \vy) &= \tvu_L \equiv (\widetilde{u}_L,0) , \ \tvu(t, L, \vy) = \tvu_R \equiv (\widetilde{u}_R,0).
\label{Eu4}
\end{align}

Our goal is to show that the planar rarefaction waves are stable in the vanishing dissipation limit for the NSF system \eqref{E1}--\eqref{E5}. Specifically, 
any family of global in time weak solutions $(\vre, \vte, \vue)_{\ep > 0}$ will approach the profile $(\tvr, \tvt, \tvu)$ as $\ep \to 0$ in the energy norm 
for any $0 < t \leq T$ as long as the initial data converge to the Riemann data \eqref{E10} in the same norm. To the best of our knowledge, this is the first result on global stability of planar rarefaction waves for the NSF system in the framework of 
global--in--time
\emph{weak} solution 
in the full 3-D setting. There are relevant 1-D results by 
Chen and Perepelitsa \cite{ChenPer1}, see also Goodman and Xin \cite{GooXin}, Hoff and Liu \cite{HofLiu} or Xin \cite{Xin93}.
Besides the afore mentioned paper by Li, Wang, and Wang \cite{LiWanWan}, the are several other results in the context of small and smooth perturbations: Li and Li 
\cite{LiLi1}, 
Li, Wang, and Wang \cite{LiHWanWan}, \cite{LiWaWa}, \cite{LiWanTWan}, 
among others.
 
The main ingredients of our approach include:
\begin{itemize}
	\item A general framework of weak solutions for the NSF system with inhomogeneous Dirichlet boundary conditions developed in \cite{FeiNovOpen} (see also \cite{ChauFei}). 
	\item Stability estimates based on the (relative) ballistic energy.
	\item Suitable choice of EOS compatible both with the Boyle--Mariotte law for the rarafeaction waves solution and the extra hypothesis required by the existence theory 
	developed in \cite{FeiNovOpen}.
	\end{itemize}

The paper is organized as follows. In Section \ref{M}, we recall the existence theory for the NSF system and state our main result. Section \ref{e} is devoted to the basic uniform 
bounds independent of the vanishing parameter $\ep$. In Section \ref{r}, we introduce the Bregman distance based on the ballistic energy functional and show convergence on compact subintervals of $(0,T]$. The proof of the main result is completed in Section \ref{P}. This last step is due to incompatibility of the ``temperature'' test function that must be smooth to obtain the energy 
estimates but coincides with the rarefaction wave in the relative energy estimates, see Remark \ref{RR1}.

\section{Basic hypotheses and main results}
\label{M}

We recall the basic structural restrictions imposed on EOS and the transport coefficients, summarize the available global existence results, and, finally, state our main result. 

\subsection{Constitutive relations}

Motivated by \eqref{E7}, we consider the pressure EOS in the form 
\begin{equation} \label{M1}
	p = p_\ep (\vr, \vt) = {\vt}^{\frac{5}{2}} P(Z) + a(\ep) \vt^4,\ Z = \frac{\vr}{\vt^{\frac{3}{2}}},\ a(\ep) > 0,
\end{equation}
The extra term $a(\ep) \vt^4$ represent the so--called radiation pressure and plays a crucial role in the existence theory, see \cite[Part II]{FeiNovOpen}. Accordingly, we have 
\begin{align} 
e = e_\ep (\vr, \vt) &= \frac{3}{2} \frac{ {\vt}^{\frac{5}{2}} }{\vr} P(Z) +   \frac{3 a(\ep)}{\vr} \vt^4   \br
s = s_\ep (\vr, \vt) &=	S(Z) + \frac{4 a(\ep)}{\vr} \vt^3, 
	\label{M2}
\end{align}
where
\begin{equation} \label{M3} 
S'(Z)  = -\frac{3}{2} \frac{ \frac{5}{3} P(Z) - P'(Z) Z}{Z^2} < 0.
\end{equation}

The existence theory developed in \cite{FeiNovOpen} requires two extra physically grounded assumptions that are, however, not compatible 
with the Boyle--Mariotte law \eqref{BM}: 
\begin{equation} \label{M4}
\frac{P(Z)}{Z^{\frac{5}{3}}} \searrow p_\infty > 0,\ S(Z) \to 0 \ \mbox{as}\ Z \to \infty.
\end{equation}
Note that \eqref{BM} affects the EOS only in the degenerate are of large parameter $Z = \frac{\vr}{\vt^{\frac{3}{2}}}$, while the fluid may still comply 
with \eqref{BM}	for moderate values of $Z$. One of the simplest choices is 
\begin{equation} \label{M5}
	P(Z) = Z \ \mbox{for}\ 0 \leq Z \leq \wtZ,\ 
		P(Z) = \frac{3}{5} Z^{\frac{5}{3}} \wtZ^{-\frac{2}{3}} + \frac{2}{5} \wtZ,\ Z > \wtZ, 
\end{equation}
\[
S'(Z) = - \frac{1}{Z} \ \mbox{for}\ Z \leq \wtZ, S'(Z) = - \frac{\wtZ}{Z^2}, \ Z > \wtZ	
\]
yielding
\begin{equation} \label{M6}
S(Z) = 1 - \log \left( \frac{Z}{\wtZ} \right) \ \mbox{if}\ 0 < Z \leq \wtZ,\ S(Z) = \frac{\wtZ}{Z} \ \mbox{if}\ Z > \wtZ.
\end{equation}
Finally, we fix $\wtZ$ related to the amplitude of the constant entropy of the rarefaction wave meaning
\begin{equation} \label{M7}
	\wtZ = \frac{\tvr}{\tvt^{\frac{3}{2}}}.
\end{equation}

As for the transport coefficients, we adopt the hypotheses of \cite[Chapter 12]{FeiNovOpen}: 
\begin{align} 
\mu \in C^1[0, \infty),\ 0 < \underline{\mu} (1 + \vt) \leq \mu(\vt) \leq \Ov{\mu}(1 + \vt), \br
\eta \in C^1[0, \infty),\ 0 \leq \eta(\vt) \leq \Ov{\eta}(1 + \vt), \ \br
\kappa \in C^1[0, \infty),\ 0 < \underline{\kappa} (1 + \vt^\beta) \leq \kappa(\vt) \leq \Ov{\kappa}(1 + \vt^\beta),\ \beta > 6.	
	\label{M8}
\end{align}

\subsection{Boundary conditions}

First, we set
\begin{align} \label{b2}
	\Gamma_{\rm in} &= \left\{ x_1 = - L \ \Big|\ 
	\vy \in \mathbb{T}^{d-1} \right\},\ \Gamma_{\rm out} = \left\{ x_1 =  L \ \Big|\ 
	\vy \in \mathbb{T}^{d-1} \right\} \ \mbox{if}\ \widetilde{u}_L > 0,
 \br
	\Gamma_{\rm in} &= \left\{ x_1 = L \ \Big|\ 
\vy \in \mathbb{T}^{d-1} \right\},\ \Gamma_{\rm out} = \left\{ x_1 = - L \ \Big|\ 
\vy \in \mathbb{T}^{d-1} \right\} \ \mbox{if}\ \widetilde{u}_R < 0.
	\end{align}
Note carefully that the rarefaction wave solution always satisfies
\begin{equation} \label{b2a}
	\partial_{x_1} \widetilde{u} \geq 0, 
	\end{equation}
see e.g. \cite{FeKr2014}.

We prescribe the following boundary conditions to be satified by the weak solutions of the NSF system:
\begin{align}
\vu(t,-L, \vy) &= \tvu_L,\ \vu(0,L,\vy) = \tvu_R, \br
\vt(t,-L, \vy) &= \tvt_L,\ \vt(0,L,\vy) = \tvt_R, \br 
\vr(t,-L, \vy) &= \tvr_L \ \mbox{if}\ \tvu_L > 0,\ 
\vr(t,L, \vy)= \tvr_R \ \mbox{if}\ \tvu_R < 0. 
	\label{BC}
	\end{align}

\subsection{Weak solutions to the NSF system}
\label{ws}

Let $\Omega$ be given by \eqref{b1}. We say that $(\vr, \vt, \vu)$ is a weak solutions of the NSF system \eqref{E1}--\eqref{E5}, 
with the boundary conditions \eqref{BC} if the following holds:
\begin{itemize}
\item {\bf Equation of continuity.} If $\tvu_L > 0$, then the integral identity 
	\begin{align}
	\int_0^\tau &\intO{ \Big[ \vr \partial_t \varphi + \vr \vu \cdot \Grad \varphi \Big]} \dt \br 
	&= - \int_0^\tau \int_{x_1 = -L} \varphi \tvr_L \tvu_L \ \D \sigma_x \dt + \int_0^\tau \int_{x_1=L} \varphi \vr \tvu_R \D \sigma_x \dt  + \left[ \intO{ \varphi \vr } \right]_{t = 0}^{t = \tau}
	\label{WS1}
\end{align}	
holds for any $0 \leq \tau \leq T$, $\varphi \in C^1([0,T] \times \Ov{\Omega})$.

If $\tvu_R < 0$, then the integral identity 
\begin{align}
	\int_0^\tau &\intO{ \Big[ \vr \partial_t \varphi + \vr \vu \cdot \Grad \varphi \Big]} \dt \br 
	&= - \int_0^\tau \int_{x_1 = -L} \varphi \vr \tvu_L \ \D \sigma_x \dt + \int_0^\tau \int_{x_1=L} \varphi \tvr_R \tvu_R \D \sigma_x \dt  + \left[ \intO{ \varphi \vr } \right]_{t = 0}^{t = \tau}
	\label{WS2}
\end{align}	
holds for any $0 \leq \tau \leq T$, $\varphi \in C^1([0,T] \times \Ov{\Omega})$.

\item {\bf Momentum equation.} 
\[
(\vu - \vuB) \in L^2(0,T; W^{1,2}_0 (\Omega; R^d)), 
\]	
where $\vuB$ is continuously differentiable function satifying 
\[
\vuB|_{x_1 = - L} = \tvu_L,\ \vuB|_{x_1 = L} = \tvu_R. 
\]
The integral identity 	
	\begin{align}
	\int_0^\tau &\intO{ \Big[ \vr \vu \cdot \partial_t \bfphi + \vr \vu \otimes \vu : \Grad \bfphi + p_\ep(\vr, \vt) \Div \bfphi \Big] } \dt \br &= 
	\int_0^\tau \intO{ \mathbb{S}_\ep (\vt, \Ds \vu) : \Ds \bfphi } \dt + \left[ \intO{ \vr \vu \cdot \bfphi } \right]_{t=0}^{t = \tau}  
	\label{WS3}
\end{align}	
for any $0 \leq \tau \leq T$, $\bfphi \in C^1_c([0,T] \times \Omega; R^d)$.

\item {\bf Entropy balance.} The integral inequality
 \begin{align}
	&\left[ \intO{ \vr s_\ep(\vr, \vt) \varphi } \right]_{t = \tau_1}^{t = \tau_2} - \int_{\tau_1}^{\tau_2} \intO{ 
		\left[ \vr s_\ep(\vr, \vt) \partial_t \varphi + \vr s_\ep (\vr, \vt) \vu \cdot \Grad \varphi + \frac{\vc{q}_\ep}{\vt} \cdot \Grad \varphi \right] } \dt \br &\geq 
	\int_{\tau_1}^{\tau_2} \intO{ \frac{\varphi}{\vt} \left( \mathbb{S}_\ep(\vt, \Ds \vu ): \Ds \vu - \frac{\vc{q}_\ep(\vt, \Grad \vt) \cdot \Grad \vt }{\vt} \right) } \dt 
	\label{WS4}
\end{align}
for any $0 \leq \tau_1 < \tau_2 \leq T$, and any $\varphi \in C^1_c([0,T] \times \Omega)$, $\varphi \geq 0$.

\item {\bf Ballistic energy inequality.} If $\tvu_L > 0$, then for any 
\[
\Theta \in C^1([0,T] \times \Ov{\Omega}),\ \Theta > 0,\ \Theta|_{x_1 = -L} = \tvt_L,\ \Theta|_{x_1 = L} = \tvt_R
\]
there holds
	\begin{align}  
	&- 		\int_0^T \partial_t \psi \intO{ \left( \frac{1}{2} \vr |\vu - \vuB|^2 + \vr e_\ep (\vr, \vt) - \Theta \vr s_\ep (\vr, \vt) \right) } \dt  \br &- 
	\int_0^T \psi \int_{x_1 = -L} \Big[ \tvr_L e_\ep (\tvr_L, \tvt_L) - \tvt_L \tvr_L s_\ep (\tvr_L, \tvt_L) \Big] \tvu_L \ \D \sigma_x \dt \br &+ 
	\int_0^T \psi \int_{x_1 = L} \Big[ \vr e_\ep (\vr, \tvt_R) - \tvt_R \vr s_\ep (\vr, \tvt_R) \Big] \tvu_R \ \D \sigma_x \dt \br &+ 
	\int_0^T \psi \intO{ \frac{\Theta}{\vt}	 \left( \mathbb{S}_\ep (\vt, \Ds \vu) : \Ds \vu - \frac{\vc{q}_\ep (\vt, \Grad \vt) \cdot \Grad \vt }{\vt} \right) } \dt \br
	&\leq	 \psi(0)  \intO{ \left( \frac{1}{2} \vr |\vu - \vuB|^2 + \vr e_\ep(\vr, \vt) - \Theta \vr s_\ep(\vr, \vt) \right) (0, \cdot) }  \br  
	&		- \int_0^T \psi \intO{ \Big[ \vr (\vu - \vuB) \otimes (\vu - \vuB) + p_\ep(\vr, \vt) \mathbb{I} - \mathbb{S}_\ep (\vt, \Ds \vu) \Big] : \Ds \vuB } \dt \br & -
	\int_0^T \psi \intO{ \vr (\vu - \vuB)\cdot (\partial_t \vuB + \vuB \cdot \Grad \vuB) } \dt \br 
	&- \int_0^T \psi \intO{ \left[ \vr s_\ep (\vr, \vt) \left( \partial_t \Theta + \vu \cdot \Grad \Theta \right) + \frac{\vc{q}_\ep(\vt, \Grad \vt) }{\vt} \cdot \Grad \Theta \right] } \dt
	\label{WS5}
\end{align}
for any $\psi \in C^1_c [0,T)$, $\psi \geq 0$.

If $\tvu_R < 0$, the inequality \eqref{WS5} is replaced by 
	\begin{align}  
	&- 		\int_0^T \partial_t \psi \intO{ \left( \frac{1}{2} \vr |\vu - \vuB|^2 + \vr e_\ep (\vr, \vt) - \Theta \vr s_\ep (\vr, \vt) \right) } \dt  \br &- 
	\int_0^T \psi \int_{x_1 = -L} \Big[ \vr e_\ep (\vr, \tvt_L) - \tvt_L \vr_L s_\ep (\vr_L, \tvt_L) \Big] \tvu_L \ \D \sigma_x \dt \br &+ 
	\int_0^T \psi \int_{x_1 = L} \Big[ \tvr_R e_\ep (\tvr_R, \tvt_R) - \tvt_R \tvr_R s_\ep (\tvr_R, \tvt_R) \Big] \tvu_R \ \D \sigma_x \dt \br &+ 
	\int_0^T \psi \intO{ \frac{\Theta}{\vt}	 \left( \mathbb{S}_\ep (\vt, \Ds \vu) : \Ds \vu - \frac{\vc{q}_\ep (\vt, \Grad \vt) \cdot \Grad \vt }{\vt} \right) } \dt \br
	&\leq	 \psi(0)  \intO{ \left( \frac{1}{2} \vr |\vu - \vuB|^2 + \vr e_\ep(\vr, \vt) - \Theta \vr s_\ep(\vr, \vt) \right) (0, \cdot) }  \br  
	&		- \int_0^T \psi \intO{ \Big[ \vr (\vu - \vuB) \otimes (\vu - \vuB) + p_\ep(\vr, \vt) \mathbb{I} - \mathbb{S}_\ep (\vt, \Ds \vu) \Big] : \Ds \vuB } \dt \br & -
	\int_0^T \psi \intO{ \vr (\vu - \vuB)\cdot (\partial_t \vuB + \vuB \cdot \Grad \vuB) } \dt \br 
	&- \int_0^T \psi \intO{ \left[ \vr s_\ep (\vr, \vt) \left( \partial_t \Theta + \vu \cdot \Grad \Theta \right) + \frac{\vc{q}_\ep(\vt, \Grad \vt)}{\vt} \cdot \Grad \Theta \right] } \dt. 
	\label{WS6}
\end{align}
	\end{itemize}

The reader may consult \cite[Chapter 12]{FeiNovOpen} for the basic properties of the weak solutions to the NSF system. In particular, we report the following existence result 
\cite[Chapter 12, Theorem 18]{FeiNovOpen}. 

\begin{Proposition}[\bf Global existence for the NSF system] \label{PP1}
Let $\Omega$ be given by \eqref{b1}. Suppose that the thermodynamic functions $p$, $e$, $s$ and the transport coefficients $\mu$, $\eta$, $\kappa$ satisfy the hypotheses \eqref{M1}--\eqref{M8}. 

Then the NSF system \eqref{E1}--\eqref{E5}, with the boundary conditions \eqref{BC}, admits a weak solution $(\vr, \vt, \vu)$ for any initial data 
\begin{align} 
\vr(0, \cdot) &= \vr_0,\ (\vr \vu)(0, \cdot) = \vr_0 \vu_0,\ \vr s(\vr, \vt)(0, \cdot) = \vr_0 s(\vr_0, \vt_0), \br
\vr_0,\ \vt_0 &\in L^\infty(\Omega),\ \vu_0 \in L^\infty(\Omega;R^d),\ \inf_\Omega \vr_0 > 0,\ \inf_\Omega \vt_0 > 0.
\label{WS7}	
	\end{align}
	\end{Proposition}

\subsection{Main result}

Having collected all necessary material, we are in a position to state our main result 

\begin{Theorem}[\bf Stability of planar rarefaction waves] \label{MT1}
	Under the hypotheses of Proposition \ref{PP1}, let $(\tvr, \tvt, \tvu)$ be a planar rarefaction wave solution of the Euler system
	introduced in Secton \ref{RWS}. Let $(\vre, \vte, \vue)_{\ep > 0}$ be a family of weak solutions of the NSF system, with the boundary conditions \eqref{BC},  and the initial data 
	\begin{align}
	0 < \underline{\vr} \leq \vr_{0,\ep} \leq \Ov{\vr},\ 0 < \underline{\vt} \leq \vt_{0,\ep} \leq \Ov{\vt}, \ 
	|\vm_{0,\ep}= \vr_{0, \ep} \vu_{0,\ep}| \leq \Ov{m}
	\label{hyp1}
	\end{align}
satisfying
\begin{equation} \label{hyp2}
	\vr_{0, \ep} \to \tvr(0, \cdot),\ \vt_{0, \ep} \to \tvt(0, \cdot),\ 
	\vr_{0, \ep} \vu_{0,\ep} \to \tvr \tvu (0, \cdot) \ \mbox{in}\ L^1(\Omega) \ \mbox{as}\ \ep \to 0.
\end{equation}	
Let 
\begin{equation} \label{hyp3}
0 < \frac{a(\ep)}{\ep} \to 0 \ \mbox{as}\ \ep \to 0.
\end{equation}
Then 
\begin{equation} \label{conc}
	\vr_{\ep}(t, \cdot) \to \tvr(t, \cdot),\ \vte(t, \cdot) \to \tvt(t, \cdot),\ 
\vre \vue (t, \cdot) \to \tvr \tvu (t, \cdot) \ \mbox{in}\ L^1(\Omega) \ \mbox{as}\ \ep \to 0 
\end{equation}
uniformly for $t$ belonging to compact sets in $(0,T]$.
	\end{Theorem}

The rest of the paper is devoted to the proof of Theorem \ref{MT1}. For definiteness, we focus on the case 
\begin{equation} \label{HHYP}
	\tvu_L > 0
\end{equation}	
seeing that the complementary case $\tvu_R < 0$ can be handled in the same way. 
	
\section{Uniform bounds}
\label{e}	

We start by deriving uniform bounds on the family $(\vre, \vte, \vue)_{\ep > 0}$ independent of the scaling parameter $\ep \to 0$.
\subsection{Mass balance}

Recalling our hypothesis \eqref{HHYP} we deduce from the weak formulation of the equation of continuity \eqref{WS1}: 
\begin{equation} \label{Mm1}
\intO{ \vre (s, \cdot) } + \int_0^z \int_{x_1 = L} \vre \tvu_R \rm \D \sigma_x \dt \leq
\intO{ \vr_{0,\ep} (0, \cdot) } - \int_0^z \int_{x_1 = - L} \vr_L \tvu_L \D \sigma_x \dt 
	\end{equation}
for any $z \in [0,T]$. This yields the following bound 
\begin{equation} \label{Mm2}
\sup_{z \in [0,T]} \left( \| \vre(z, \cdot)  \|_{L^1(\Omega)} + \| \vre(z, \cdot) \|_{L^1(x_1 = L)} \right) \aleq 1.
\end{equation}
The reader may consult \cite[Chapter 3, Section 3.3.1]{FeiNovOpen} for the interpretation of the ``trace'' of the density $\vre$ on the outflow part of the domain.

\subsection{Energy estimates}

Our next goal is to derive uniform bounds from the ballistic energy inequality \eqref{WS5}. To this end, we fix 
$\vuB$ satisfying the boundary conditions
\begin{equation} \label{e2}
\vuB = [u_B(x_1), 0],\ \partial_{x_1} u_B \geq 0, \ u_B(-L) = \tvu_L,\ u_B(L) = \tvu_R, 
\end{equation}
and $\Theta = \vtB$ -- a superharmonic function
\begin{equation} \label{e3} 
\Del \vtB \geq 0 \ \mbox{in}\ \Omega,\ \vtB(-L) = \tvt_L,\ \vtB(L) = \tvt_R.
\end{equation}

\begin{Remark} \label{RR1}
A proper choice of $\vtB$ is absolutely crucial in the proof. As the boundary values $\tvt_R$, $\tvt_L$ are positive constants and the travelling wave component 
$\tvt$ strictly positive, we may consider
\[
\vtB = \vtB(x_1), \ \vtB(-L) = \tvt_L,\ \vtB(L) = \tvt_R,\ 0 < \vtB \leq \tvt,\ \partial^2_{x_1,x_1} \vtB \geq 0.
\] 	
	\end{Remark}

Plugging the above anstaz in the ballistic energy inequality \eqref{WS5} we obtain:
	\begin{align}
	\Big[ & \intO{ \left( \frac{1}{2} \vre |\vue - \vuB |^2 + \vre e_\ep (\vre, \vte) - \vtB \vre s_\ep (\vre, \vte) \right) } \Big]_{t=0}^{t=z} 
	\br 
	&+ \int_0^z \int_{x_1 = L} \Big( \vre e_\ep (\vre, \tvt_R) - \tvt_R \vre s(\vre, \tvt_R) \Big) \tvu_R \ \D \sigma_x \dt \br 
	& - \int_0^z \int_{x_1 = -L}
	\Big( \tvr_L e_\ep (\tvr_L, \tvt_L) - \tvt_L \tvr_L s_\ep (\tvr_L, \tvt_L) \Big) \tvu_L \ \D \sigma_x \dt  \br
	&+ \int_0^z \intO{ \frac{\vtB}{\vte}  \left( \mathbb{S}_\ep(\vte, \Ds \vue) : \Ds \vue + \ep \frac{\kappa(\vte) |\Grad \vte|^2 }{\vte} \right) } \dt \br
	&\leq - \int_0^z \intO{ \Big[ \vre (\vue -  \vuB) \otimes (\vue - \vuB) + p_\ep(\vre, \vte) \mathbb{I} - \mathbb{S}_\ep (\vte, \Ds \vue) \Big] : \Ds \vuB } \dt \br  &\quad - \int_0^z \intO{ \vre \Big[ (\vuB \cdot \Grad) \vuB \Big] \cdot (\vue - \vuB) } \dt \br 
	&\quad - \int_0^z \intO{ \Big[ \vre s_\ep (\vre, \vte) \vue \cdot \Grad \vtB - \ep \frac{\kappa(\vte) \Grad \vte }{\vte} \cdot \Grad \vtB \Big] } \dt .
	\label{e1}
\end{align}

First observe that \eqref{e2} yields 
\[
- \int_0^z \intO{ \Big[ \vre (\vue -  \vuB) \otimes (\vue - \vuB) + p_\ep(\vre, \vte) \mathbb{I} \Big] : \Ds \vuB } \dt \leq 0.
\]

Next, by virtue of hypothesis \eqref{M8} and Korn--Poincar\' e inequality, 
\begin{equation} \label{ee1}
\ep \| \vue \|^2_{W^{1,2}(\Omega; R^d)} \aleq \left( \ep  + \intO{ \frac{\vtB}{\vte}  \mathbb{S}_\ep(\vte, \Ds \vue) : \Ds \vue  } \right). 
\end{equation}
Similarly, since $\vt$ admits a bounded trace on $\partial \Omega$, we may use Poincar\' e inequality together with hypothesis \eqref{M8} to obtain 
\begin{equation} \label{e6}
	\| \log (\vte) \|_{L^2(\Omega)}^2 +	\| \vte^{\frac{\beta}{2}} \|_{W^{1,2}(\Omega)}^2 \aleq 	\intO{ \frac{\kappa (\vte)}{\vte^2}|\Grad \vte|^2 } + 1,\ \beta > 6.
\end{equation}

Consequently, inequality \eqref{e1} gives rise to 
	\begin{align}
	\Big[ & \intO{ \left( \frac{1}{2} \vre |\vue - \vuB |^2 + \vre e_\ep (\vre, \vte) - \vtB \vre s_\ep (\vre, \vte) \right) } \Big]_{t=0}^{t=z} 
	\br 
	&+ \int_0^z \int_{x_1 = L} \Big( \vre e_\ep (\vre, \tvt_R) - \tvt_R \vre s(\vre, \tvt_R) \Big) \tvu_R \ \D \sigma_x \dt \br 
	&+ \ep \int_0^z \| \vue \|^2_{W^{1,2}(\Omega; R^d)} \dt +   \ep \int_0^z \left( \| \log (\vte) \|_{L^2(\Omega)}^2 +	\| \vte^{\frac{\beta}{2}} \|_{W^{1,2}(\Omega)}^2 \right)  \dt \br
	&\aleq \int_0^z \intO{ \mathbb{S}_\ep (\vte, \Ds \vue) : \Ds \vuB } \dt  - \int_0^z \intO{ \vre \Big[ (\vuB \cdot \Grad) \vuB \Big] \cdot (\vue - \vuB) } \dt \br 
	&\quad - \int_0^z \intO{ \Big[ \vre s_\ep (\vre, \vte) \vue \cdot \Grad \vtB - \ep \frac{\kappa(\vte) \Grad \vte }{\vte} \cdot \Grad \vtB \Big] } \dt  + z. 
	\label{e1a}
\end{align}

\subsubsection{Integrals containing velocity}

By virtue of hypothesis \eqref{M8}, 
\[
\intO{ \mathbb{S}_\ep (\vte, \Ds \vue) : \Ds \vuB } \leq \delta \ep \| \vue \|^2_{W^{1,2}(\Omega; R^d)} + c(\delta) \ep \intO{ \vte^2 }.
\]
Similarly, 
\[
\intO{ \vre \Big[ (\vuB \cdot \Grad) \vuB \Big] \cdot (\vue - \vuB) } \leq \delta \intO{ \vre |\vue - \vuB|^2 } + c(\delta) \intO{ \vre },
\]
where $\delta > 0$ is arbitrary. Consequently, using the uniform bounds \eqref{Mm2} we may rewrite \eqref{e1a} in the following form 
	\begin{align}
	\Big[ & \intO{ \left( \frac{1}{2} \vre |\vue - \vuB |^2 + \vre e_\ep (\vre, \vte) - \vtB \vre s_\ep (\vre, \vte) \right) } \Big]_{t=0}^{t=z} 
	\br 
	&+ \int_0^z \int_{x_1 = L} \Big( \vre e_\ep (\vre, \tvt_R) - \tvt_R \vre s(\vre, \tvt_R) \Big) \tvu_R \ \D \sigma_x \dt \br 
	&+ \frac{\ep}{2} \int_0^z \| \vue \|^2_{W^{1,2}(\Omega; R^d)} \dt +   \frac{\ep}{2} \int_0^z \left( \| \log (\vte) \|_{L^2(\Omega)}^2 +	\| \vte^{\frac{\beta}{2}} \|_{W^{1,2}(\Omega)}^2 \right)  \dt  \br
	&\aleq 
	- \int_0^z \intO{ \Big[ \vre s_\ep (\vre, \vte) \vue \cdot \Grad \vtB - \ep \frac{\kappa(\vte) \Grad \vte }{\vte} \cdot \Grad \vtB \Big] } \dt  \br 
	&\quad + \int_0^z \intO{ \vre |\vue - \vuB|^2  } \dt  + z. 
	\label{e1b}
\end{align}

\subsubsection{Integrals containing temperature}

First, we use the fact that $\vtB$ is superharmonic and integrate by parts
\begin{align} 
&\intO{ \frac{\kappa(\vte) \Grad \vte }{\vte} \cdot \Grad \vtB } = - \int_{x_1 = -L} K(\tvt_L) \partial_{x_1} \vtB\  \D \sigma_x  + 
\int_{x_1 = L} K(\tvt_R) \partial_{x_1} \vtB \ \D \sigma_x \br  
&- \intO{ K(\vte) \Del \vtB } \leq
- \int_{x_1 = -L} K(\tvt_L) \partial_{x_1} \vtB\  \D \sigma_x  + 
\int_{x_1 = L} K(\tvt_R) \partial_{x_1} \vtB \ \D \sigma_x,\
 \ K'(Y) = \frac{\kappa(Y)}{Y}.
\label{e1c}
\end{align}

Next, in accordance with hypothesis \eqref{M2}, \eqref{M6}, we write 
\[
\intO{ \vre s_\ep (\vre, \vte) \vue \cdot \Grad \vtB } = \intO{ \vre S \left( \frac{\vre}{\vte^{\frac{3}{2}}} \right) \vue \cdot \Grad \vtB } + 
4 \ep \intO{ \vte^3 \vue \cdot \Grad \vtB }.
\]
Using \eqref{e6} we obtain 
\[
\intO{ \vte^3 \vue \cdot \Grad \vtB } \leq \delta \| \vue \|^2_{L^2(\Omega; R^d)} + c(\delta) \intO{ \vte^6 }
\]
for any $\delta > 0$. Since $\beta > 6$, we may use \eqref{e6} to absorb this term by the integral on the left--hand side of \eqref{e1b}. 
Accordingly, relation \eqref{e1b} reduces to 
	\begin{align}
	\Big[ & \intO{ \left( \frac{1}{2} \vre |\vue - \vuB |^2 + \vre e_\ep (\vre, \vte) - \vtB \vre s_\ep (\vre, \vte) \right) } \Big]_{t=0}^{t=z} 
	\br 
	&+ \int_0^z \int_{x_1 = L} \Big( \vre e_\ep (\vre, \tvt_R) - \tvt_R \vre s(\vre, \tvt_R) \Big) \tvu_R \ \D \sigma_x \dt \br 
	&+ \frac{\ep}{2} \int_0^z \| \vue \|^2_{W^{1,2}(\Omega; R^d)} \dt +   \frac{\ep}{2} \int_0^z \left( \| \log (\vte) \|_{L^2(\Omega)}^2 +	\| \vte^{\frac{\beta}{2}} \|_{W^{1,2}(\Omega)}^2 \right)  \dt \br
	&\aleq 
	- \int_0^z \intO{ \vre S \left( \frac{\vre}{\vte^{\frac{3}{2}}} \right) \vue \cdot \Grad \vtB } \dt  + \int_0^z \intO{ \vre |\vue - \vuB|^2  } \dt  + z.
	\label{e1d}
\end{align}

Finally, it remains to control the first integral on the right--hand side of \eqref{e1d}. To this end, we make use of hypothesis \eqref{M6}. 
First, suppose 
\[
	\frac{\vre}{\vte^{\frac{3}{2}}} > \wtZ.
\]	
It follows from \eqref{M6} that
\begin{equation} \label{e13}
\left| \vre S \left( \frac{\vre}{\vte^{\frac{3}{2}}} \right) \vue \cdot \Grad \vtB \right| \aleq  \vre |\vue - \vuB|^2 + \vre. 
\end{equation}
If 
\[
\frac{\vre}{\vte^{\frac{3}{2}}} \leq \wtZ,
\]	
then, by the same token 
\begin{align}
\left| \vre S \left( \frac{\vre}{\vte^{\frac{3}{2}}} \right) \vue \cdot \Grad \vtB \right| &\aleq 
\vre |\log(\vre)| |\vue| + \vre |\log^+(\vte)| |\vue| \br &\aleq 
\vre + \vre |\vue - \vuB|^2 + \vre |\log(\vre)|^2 + \vre |\log^+(\vte)|^2.
\label{e14}
\end{align}
Summing up \eqref{e13}, \eqref{e14} we rewrite \eqref{e1d} in the form 
	\begin{align}
	\Big[ & \intO{ \left( \frac{1}{2} \vre |\vue - \vuB |^2 + \vre e_\ep (\vre, \vte) - \vtB \vre s_\ep (\vre, \vte) \right) } \Big]_{t=0}^{t=z} 
	\br 
	&+ \int_0^z \int_{x_1 = L} \Big( \vre e_\ep (\vre, \tvt_R) - \tvt_R \vre s(\vre, \tvt_R) \Big) \tvu_R \ \D \sigma_x \dt \br 
	&+ \frac{\ep}{2} \int_0^z \| \vue \|^2_{W^{1,2}(\Omega; R^d)} \dt +   \frac{\ep}{2} \int_0^z \left( \| \log (\vte) \|_{L^2(\Omega)}^2 +	\| \vte^{\frac{\beta}{2}} \|_{W^{1,2}(\Omega)}^2 \right)  \dt \br
	&\aleq 
	 \int_0^z \intO{ \left( \vre |\vue - \vuB|^2 +  \vre^{\frac{5}{3}} + \vre \vte \right) } \dt  + z.
	\label{e1e}
\end{align}
Now, we exploit hypotheses \eqref{M2}, \eqref{M5} to obtain the final estimate 
	\begin{align}
	\Big[ & \intO{ \left( \frac{1}{2} \vre |\vue - \vuB |^2 + \vre e_\ep (\vre, \vte) - \vtB \vre s_\ep (\vre, \vte) \right) } \Big]_{t=0}^{t=z} 
	\br 
	&+ \int_0^z \int_{x_1 = L} \Big( \vre e_\ep (\vre, \tvt_R) - \tvt_R \vre s(\vre, \tvt_R) \Big) \tvu_R \ \D \sigma_x \dt \br 
	&+ \frac{\ep}{2} \int_0^z \| \vue \|^2_{W^{1,2}(\Omega; R^d)} \dt +   \frac{\ep}{2} \int_0^z \left( \| \log (\vte) \|_{L^2(\Omega)}^2 +	\| \vte^{\frac{\beta}{2}} \|_{W^{1,2}(\Omega)}^2 \right)  \dt \br
	&\aleq 
	\int_0^z \intO{ \left( \frac{1}{2} \vre |\vue - \vuB|^2 + \vre e_\ep (\vre, \vte) - \vtB \vre s_\ep (\vre, \vte) \right) } \dt  + z.
	\label{e1f}
\end{align}
Note we have used the fact the entropy is dominated by the energy. Indeed, similarly to \eqref{e14}, 
\begin{align} 
	0 \leq \vr s_\ep (\vr, \vt) &\aleq a(\ep) \vt^3 + \vr (1 + |\log(\vr)|) + \vr \log^+ (\vt),\br 
	\vr e_\ep (\vr, \vt) &\ageq a(\ep) \vt^4 + \vr^{\frac{5}{3}} + \vr \vt.
	\label{e16} 
\end{align}

Thus we may apply Gronwall's argument to \eqref{e1f}, which, combined with the estimate \eqref{Mm2}, yields the following result.

\begin{Lemma}[\bf Uniform bounds] \label{eL1}
	Under the hypotheses of Theorem \ref{MT1}, 
\begin{align}
	&\intO{ \left( \frac{1}{2} \vre |\vue - \vuB |^2 + \vre e_\ep (\vre, \vte) - \vtB \vre s_\ep (\vre, \vte) \right)(z, \cdot) } \br
	&+ \frac{\ep}{2} \int_0^z \| \vue \|^2_{W^{1,2}(\Omega; R^d)} \dt +   \frac{\ep}{2} \int_0^z \left( \| \log (\vte) \|_{L^2(\Omega)}^2 +	\| \vte^{\frac{\beta}{2}} \|_{W^{1,2}(\Omega)}^2 \right)  \dt \br 
	&\leq \intO{ \left( \frac{1}{2} \vr_{0, \ep} |\vu_{0,\ep} - \vuB |^2 + \vr_{0,\ep} e_\ep (\vr_{0, \ep}, \vt_{0,\ep}) - \vtB \vr_{0, \ep} s_\ep (\vr_{0, \ep}, \vt_{0, \ep} ) \right) } + Cz
	\label{e17} 
	\end{align}	
for any $0 \leq z \leq T$, the functions $\vuB$, $\vtB$ satisfy \eqref{e2}, \eqref{e3}, and the constant $C$ depends solely on the data, in particular, it is independent of $\ep > 0$.	
	\end{Lemma}

\section{Relative energy estimates}
\label{r}

Our next task is to estimate the distance between $(\vre, \vte, \vu)$ and the limit rarefaction wave profile $(\tvr, \tvt, \tvu)$.

\subsection{Relative energy}

A suitable quantity to measure the distance of a weak solutions $(\vr, \vt, \vu)$ to any trio $(\tvr, \tvt, \tvu)$ is the 
\emph{relative energy}
\begin{align}
	E &\left( \vr, \vt, \vu \Big| \tvr , \tvt, \tvu \right) \br &= \frac{1}{2}\vr |\vu - \tvu|^2 + \vr e(\vr, \vt) - \tvt \Big(\vr s(\vr, \vt) - \tvr s(\tvr, \tvt) \Big) \br&- 
	\Big( e(\tvr, \tvt) - \tvt s(\tvr, \tvt) + \frac{p(\tvr, \tvt)}{\tvr} \Big)
	(\vr - \tvr) - \tvr e (\tvr, \tvt) \br &= 
	\frac{1}{2}\vr |\vu - \tvu|^2 + \vr e(\vr, \vt) - \tvt \vr s (\vr, \vt) - \Big( e(\tvr, \tvt) - \tvt s(\tvr, \tvt) + \frac{p(\tvr, \tvt)}{\tvr} \Big) \vr + p(\tvr, \tvt).
	\label{r1}
\end{align}
As shown in \cite[Chapter 3, Section 3.1]{FeiNovOpen}, the relative energy represents a Bregman distance between $(\vr, \vt, \vu)$ and $(\tvr, \tvt, \tvu)$, 
specifically, 
\begin{itemize}
	\item 
	\[
	E \left( \vr, \vt, \vu \Big| \tvr , \tvt, \tvu \right) \geq 0;
	\]
	\item if $\vr > 0$, then 
	\[
	E \left( \vr, \vt, \vu \Big| \tvr , \tvt, \tvu \right) = 0 \ \Leftrightarrow \ 
	(\vr, \vt, \vu) = (\tvr, \tvt, \tvu);
	\]
	\item 
	Given $\tvr , \tvt, \tvu$, $E$ is a strictly convex function of the conservative entropy variables $\vr$, $\vm = \vr \vu$, $\mathcal{S} = \vr s(\vr, \vt)$.
	\end{itemize}

\subsection{Relative energy inequality}

First, we introduce a modified relative energy taking into account the radiative components of thermodynamics functions:
\[
E_\ep \left( \vr, \vt, \vu \Big| \tvr , \tvt, \tvu \right) = 
\frac{1}{2}\vr |\vu - \tvu|^2 + \vr e_\ep(\vr, \vt) - \tvt \vr s_\ep (\vr, \vt) - \Big( e_\ep (\tvr, \tvt) - \tvt s_\ep (\tvr, \tvt) + \frac{p_\ep (\tvr, \tvt)}{\tvr} \Big) \vr + p_\ep(\tvr, \tvt). 
\]

As shown in \cite[Chapter 12, Section 12.3.2]{FeiNovOpen}, any weak solution of the NSF system 
satisfies a relative energy inequality. Recalling $\tvu_L > 0$ and the fact that the rarefaction wave solution is Lipschitz on any interval $[z,T]$, 
$0 < z \leq T$, the relative energy inequality takes the form:
\begin{align}
	&\left[ \intO{ E_\ep \left(\vre, \vte, \vue \Big| \tvr, \tvt, \tvu \right) } \right]_{t = z}^{t = \tau} \br 
	&- \int_z^\tau \int_{x_1 = -L} \Big( \tvr_L e_\ep(\tvr_L, \tvt_L) - \tvt_L \tvr_L s_\ep (\tvr_L, \tvt_L) \Big) \tvu_L \ \D \sigma_x \dt \br &+ \int_z^\tau \int_{x_1 = -L} \left( e_\ep (\tvr_L, \tvt_L) - \tvt_L s_\ep (\tvt_L, \tvt_B) + \frac{p_\ep (\tvr_L, \tvt_L)}{\tvr_L} \right) \tvr_L  \tvu_L \D \sigma_x \dt
	\br 
	&+ \int_z^\tau \int_{x_1 = L} \Big( \vre e_\ep (\vre, \tvt_R) - \tvt_R \vre s_\ep (\vre, \tvt_R) \Big) \tvu_R \ \D \sigma_x \dt \br &- \int_z^\tau \int_{x_1 = L} \left( e_\ep(\tvr_R, \tvt_R) - \tvt_R s_\ep(\tvr_R, \tvt_R) + \frac{p_\ep(\tvr_R, \tvt_R)}{\tvr_R} \right) \vre \tvu_R \ \D \sigma_x \dt \br 
	&+ \int_z^\tau \intO{ \frac{\tvt}{\vte} \left( \mathbb{S}_\ep (\vte, \Ds \vue) : \Ds \vue + \ep \frac{\kappa(\vte) |\Grad \vte|^2}{\vte} \right) } \dt \br 
	&\leq  \int_z^\tau \intO{ \frac{\vre}{\tvr} (\vue - \tvu ) \cdot \Grad p_\ep(\tvr, \tvt) } \dt \br 
	&- \int_z^\tau \intO{ \left( \vre (s_\ep(\vre, \vte) - s_\ep(\tvr, \tvt)) \partial_t \tvt + \vre (s_\ep (\vre, \vte) - s_\ep(\tvr, \tvt)) \vue \cdot \Grad \tvt \right) } \dt  \br &+ \ep \int_z^\tau \intO{  
		\left( \frac{\kappa(\vte) \Grad \vte}{\vte} \right) \cdot \Grad \tvt  } \dt \br 
	&- \int_z^\tau \intO{ \Big[ \vre (\vue - \tvu) \otimes (\vue - \tvu) + p_\ep(\vre, \vte) \mathbb{I} - \mathbb{S}_\ep(\vte, \Ds \vue) \Big] : \Ds \tvu } \dt \br 
	&+ \int_z^\tau \intO{ \vre \left[  \partial_t \tvu + (\tvu \cdot \Grad) \tvu + \frac{1}{\tvr} \Grad p_\ep(\tvr, \tvt) \right] \cdot (\tvu - \vue) } \dt \br 
	&+ \int_z^\tau \intO{ \left[ \left( 1 - \frac{\vre}{\tvr} \right) \partial_t p_\ep(\tvr, \tvt) - \frac{\vre}{\tvr} \vue \cdot \Grad p_\ep(\tvr, \tvt) \right] } \dt
	\label{r2}
\end{align}

\noindent
for a.a. $\tau \geq z$.

Now observe that, in view of the uniform bounds established in Lemma \ref{eL1} and hypothesis \eqref{hyp3}, we may replace $p_\ep$, $e_\ep$, and $s_\ep$ by 
$p$, $e$, $s$ modulo an error vanishing for $\ep \to 0$. In addition, as $(\tvr, \tvu, \tvs)$ is a (Lipschitz) solution of the Euler system, we have 
\[
\partial_t \tvu + (\tvu \cdot \Grad) \tvu + \frac{1}{\tvr} \Grad p(\tvr, \tvt) = 0.
\]	 
Consequently, we may rewrite \eqref{r2} in the form
\begin{align}
	&\left[ \intO{ E_\ep \left(\vre, \vte, \vue \Big| \tvr, \tvt, \tvu \right) } \right]_{t = z}^{t = \tau} \br 
	&- \int_z^\tau \int_{x_1 = -L} \Big( \tvr_L e(\tvr_L, \tvt_L) - \tvt_L \vrB s (\tvr_L, \tvt_L) \Big) \tvu_L \ \D \sigma_x \dt \br &+ \int_z^\tau \int_{x_1 = -L} \left( e (\tvr_L, \tvt_L) - \tvt_L s (\tvt_L, \tvt_B) + \frac{p (\tvr_L, \tvt_L)}{\tvt_L} \right) \tvt_L  \tvu_L \D \sigma_x \dt
	\br 
	&+ \int_z^\tau \int_{x_1 = L} \Big( \vre e (\vre, \tvt_R) - \tvt_R \vre s (\vre, \tvt_R) \Big) \tvu_R \ \D \sigma_x \dt \br &- \int_z^\tau \int_{x_1 = L} \left( e(\tvr_R, \tvt_R) - \tvt_R s(\tvr_R, \tvt_R) + \frac{p(\tvr_R, \tvt_R)}{\tvr_R} \right) \vre \tvu_R \ \D \sigma_x \dt \br 
	&+ \int_z^\tau \intO{ \frac{\tvt}{\vte} \left( \mathbb{S}_\ep (\vte, \Ds \vue) : \Ds \vue + \ep \frac{\kappa(\vte) |\Grad \vte|^2}{\vte} \right) } \dt \br 
	&\leq  \int_z^\tau \intO{ \frac{\vre}{\tvr} (\vue - \tvu ) \cdot \Grad p(\tvr, \tvt) } \dt \br 
	&- \int_z^\tau \intO{ \left( \vre (s(\vre, \vte) - s(\tvr, \tvt)) \partial_t \tvt + \vre (s (\vre, \vte) - s(\tvr, \tvt)) \vue \cdot \Grad \tvt \right) } \dt  \br &+ \ep \int_z^\tau \intO{  
		\left( \frac{\kappa(\vte) \Grad \vte}{\vte} \right) \cdot \Grad \tvt  } \dt \br 
	&- \int_z^\tau \intO{ \Big[ \vre (\vue - \tvu) \otimes (\vue - \tvu) + p(\vre, \vte) \mathbb{I} - \mathbb{S}_\ep(\vte, \Ds \vue) \Big] : \Ds \tvu } \dt \br 
	&+ \int_z^\tau \intO{ \left[ \left( 1 - \frac{\vre}{\tvr} \right) \partial_t p(\tvr, \tvt) - \frac{\vre}{\tvr} \vue \cdot \Grad p(\tvr, \tvt) \right] } \dt
	+ \omega(\ep,z), 
	\label{r3}
\end{align}
where 
\begin{equation} \label{r3a}
	\omega(\ep, z) \to 0 \ \mbox{as}\ \ep \to 0 \ \mbox{uniformly for}\ z \ \mbox{in compact subintervals of}\ (0,T].
\end{equation}
\subsection{Integrals containing the thermodynamics functions}

We focus on the integrals
\begin{align}
	I&= \int_z^\tau \intO{ \frac{\vre}{\tvr} (\vue - \tvu ) \cdot \Grad p(\tvr, \tvt) } \dt \br 
	&- \int_z^\tau \intO{ \left( \vre (s(\vre, \vte) - s(\tvr, \tvt)) \partial_t \tvt + \vr (s(\vre, \vte) - s(\tvr, \tvt)) \vue \cdot \Grad \tvt \right) } \dt  \br&- \int_z^\tau \intO{ \Big[ \vre (\vue - \tvu) \otimes (\vue - \tvu) + p(\vre, \vt) \mathbb{I}  \Big] : \Ds \tvu } \dt \br
	&+ \int_z^\tau \intO{ \left[ \left( 1 - \frac{\vre}{\tvr} \right) \partial_t p(\tvr, \tvt) - \frac{\vre}{\tvr} \vu \cdot \Grad p(\tvr, \tvt) \right] } \dt \br
&=  - \int_z^\tau \intO{ \left( \vr (s(\vre, \vte) - s(\tvr, \tvt)) \partial_t \tvt + \vr (s(\vre, \vte) - s(\tvr, \tvt)) \vue \cdot \Grad \tvt \right) } \dt  \br &\quad - \int_z^\tau \intO{ \Big[ \vre (\vue - \tvu) \otimes (\vue - \tvu) + p(\vre, \vt) \mathbb{I}  \Big] : \Ds \tvu } \dt \br
	 	&\quad + \int_z^\tau \intO{ \left[ \left( 1 - \frac{\vre}{\tvr} \right) \partial_t p(\tvr, \tvt) - \frac{\vre}{\tvr} \tvu \cdot \Grad p(\tvr, \tvt) \right] } \dt.
	\label{r4}
\end{align}
This can be rewritten in the form 
\begin{align}
	I&=  
	- \int_z^\tau \intO{ \left( \vre (s (\vre, \vte) - s(\tvr, \tvt)) \partial_t \tvt + \vre (s(\vre, \vte) - s(\tvr, \tvt)) \tvu \cdot \Grad \tvt \right) } \dt  \br&- \int_z^\tau \intO{  \vre (\vue - \tvu) \otimes (\vue - \tvu) : \Ds \tvu  } \dt \br
	&+ \int_z^\tau \intO{ \left( 1 - \frac{\vre}{\tvr} \right) \Big( \partial_t p(\tvr, \tvt) + \tvu \cdot \Grad p(\tvr, \tvt) \Big) } \dt, \br 
	&+ \int_z^\tau \intO{ \Big(p(\tvr, \tvt) - p(\vr, \vt) \Big) \Div \tvu } - \int_{s}^\tau \int_{\partial \Omega} p(\tvr, \tvt) 
	\tvu \cdot \vc{n} \ \D \sigma_x \br 
	&+ \int_z^\tau \intO{ \vre (s(\vre, \vte) - s(\tvr, \tvt)) (\vue - \tvu) \cdot \Grad \tvt }.
	\label{r5}
\end{align}

Now, since $(\tvr, \tvt, \tvu)$ solve the Euler system, we get 
\begin{equation} \label{r6}
\partial_t \tvt + \tvu \cdot \Grad \tvt = - \frac{1}{\tvr c_v(\tvr, \tvt)} 
\tvt \partial_\vt p(\tvr, \tvt) \Div \tvu, \ \mbox{where}
\ c_v (\tvr, \tvt) = \partial_\vt e(\tvr, \tvt). 
\end{equation}
Moreover, in accordance with hypothesis \eqref{M5},
\begin{equation} \label{r7}
	p = \frac{2}{3} \vr e \ \Rightarrow\ 
\partial_\vt p(\tvr, \tvt) = \frac{2}{3} \tvr 	\partial_\vt e(\tvr, \tvt) = \frac{2}{3} \tvr c_v(\tvr, \tvt).
	\end{equation}
Consequently, applying \eqref{M5} once more we may infer that 	 
\begin{equation} \label{r8}
	\partial_t \tvt + \tvu \cdot \Grad \tvt = - \frac{2}{3} \tvt \Div \tvu. 
\end{equation}
Finally,
\begin{align}
	\partial_t &p(\tvr, \tvt) + \tvu \cdot \Grad p(\tvr, \tvt)= 
	\frac{2}{3} \Big( \partial_t (\tvr e (\tvr, \tvt)) + \tvu \cdot \Grad (\tvr e(\tvr, \tvt)) \Big) \br
	&= 	\frac{2}{3} \Big( \partial_t (\tvr e (\tvr, \tvt)) + \Div ( \tvr e(\tvr, \tvt) \tvu ) - \tvr e(\tvr, \tvt) \Div \tvu ) \Big) 
	= - \frac{2}{3} \Big( p(\tvr, \tvt) + \tvr e(\tvr, \tvt) \Big) \Div \tvu \br 
	&= - \frac{5}{3} p(\tvr, \tvt) \Div \tvu.
	\nonumber
	\end{align}
Summarizing the previous discussion, we can write \eqref{r5}
in the form	
\begin{align}
	I&=  
	\frac{2}{3} \int_z^\tau \intO{ \vre (s(\vre, \vte) - s(\tvr, \tvt))  \tvt \Div\tvu } \dt  \br&- \int_z^\tau \intO{  \vre (\vue - \tvu) \otimes (\vue - \tvu) : \Ds \tvu  } \dt \br
	&+ \frac{5}{3} \int_z^\tau \intO{ \left( \frac{\vre}{\tvr} - 1 \right) p(\tvr, \tvt) \Div \tvu  } \dt, \br 
	&+ \int_z^\tau \intO{ \Big(p(\tvr, \tvt) - p(\vre, \vte) \Big) \Div \tvu } - \int_{z}^\tau \int_{\partial \Omega} p(\tvr, \tvt) 
	\tvu \cdot \vc{n} \ \D \sigma_x \br 
	&+ \int_z^\tau \intO{ \vre (s (\vre, \vte) - s(\tvr, \tvt)) (\vue - \tvu) \cdot \Grad \tvt }.
	\label{r9}
\end{align}	

At this stage, we recall the basic properties of the rarefaction wave solutions (see e.g. \cite[formulae (4.2), (4.3)]{FeKrVa}, namely, 
\begin{equation} \label{r12}
\tvu = (\widetilde{u},0,0),\ \Div \tvu = \partial_{x_1} \widetilde{u} \geq 0,\ 
\tvu \ne 0 \ \Rightarrow \ 	\left| \frac{\partial_{x_1}\tvt}{\partial_{x_1} \widetilde{u}} \right|^2 \leq 
\frac{4}{15} \tvt.
\end{equation}

Going back to \eqref{r9}
consider the expression 
\begin{align}
D &= \frac{2}{3} \vre (s(\vre, \vte) - s(\tvr, \tvt)) \tvt \partial_{x_1} \widetilde{u} - 
\vre |u_{1,\ep} - \widetilde{u}|^2 \partial_{x_1} \widetilde{u} + \frac{5}{3} \vre \frac{p(\tvr, \tvt)}{\tvr} \partial_{x_1} \widetilde{u} \br &- \frac{2}{3} p(\tvr, \tvt) \partial_{x_1} \widetilde{u} - p(\vre, \vte) \partial_{x_1} \widetilde{u} + 
\vre (s(\vre, \vte) - s(\tvr, \tvt))(u_{1,\ep} - \widetilde{u}) \partial_{x_1}\tvt \br 
&= \frac{2}{3}\vre \partial_{x_1} \widetilde{u} \left( 
 \tvt \Big(s (\vre, \vte) - s(\tvr, \tvu) \Big) - \frac{3}{2}|u_1 - \widetilde{u}|^2 + 
\frac{5}{3}  e(\tvr, \tvt) \right) \br 
&+ \frac{2}{3} \vre \partial_{x_1} \widetilde{u} \left( - \frac{2}{3} \frac{\tvr}{\vre} e(\tvr, \tvt) -  e(\vre, \vte)   +  \frac{3}{2}
(s(\vre, \vte) - s(\tvr, \tvt))(u_{1, \ep} - \widetilde{u}) \frac{\partial_{x_1}\tvt}{\partial_{x_1} \widetilde{u}} \right)
	\label{r10}
	\end{align}
As we have no control on the derivatives of $(\tvr, \tvt, \tvu)$ for $z \to 0$, our goal is to show that $D \leq 0$. 	
	
First observe that
\[
\frac{3}{2}
(s(\vr, \vt) - s(\tvr, \tvt))(u_{1,\ep} - \widetilde{u}) \frac{\partial_{x_1}\tvt}{\partial_{x_1} \widetilde{u}} \leq \frac{3}{2} |u_{1,\ep} - \widetilde{u}|^2 + \frac{3}{8} |s(\vr, \vt) - s(\tvr, \tvt) |^2 \left| \frac{\partial_{x_1}\tvt}{\partial_{x_1} \widetilde{u}} \right|^2.
\]
Consequently, due to \eqref{r12} it is enough to show that
\begin{align}
	\tvt \Big(s(\vr, \vt) - s(\tvr, \tvt) \Big)  + 
	\frac{5}{3}  e(\tvr, \tvt)  - \frac{2}{3} \frac{\tvr}{\vr} e(\tvr, \tvt) -  e(\vr, \vt) + \frac{1}{10} |s(\vr, \vt) - s(\tvr, \tvt) |^2 \tvt \leq 0. 
	\label{r13}
\end{align}
In view of the specific form of EOS given by \eqref{M5}, \eqref{M6}, it is convenient to rewrite \eqref{r13} in terms of the variables 
$\vr$, $Z = \frac{\vr}{\vt^{\frac{3}{2}}}$. In addition, as the rarefaction wave solution satisfies 
\[
\frac{\tvr}{\tvt^{\frac{3}{2}}} = \wtZ,\ \mbox{and}\ e(\tvr, \tvt) = \frac{3}{2} \tvt,  
\] 
the inequality \eqref{r13} reduces to showing 
\begin{equation} \label{r14}
S(Z) - S(\wtZ) + \frac{5}{2} - \frac{\tvr}{\vr}	- \frac{3}{2} \frac{P(Z)}{Z^{\frac{5}{3}}} \left( \frac{\vr}{\tvr} \right)^{\frac{2}{3}} \wtZ^{\frac{2}{3}} + 
\frac{1}{10} |S(Z) - S(\wtZ) |^2 \leq 0
	\end{equation}
for any $Z > 0$, $\vr > 0$.

As $P(Z) = Z$ for $Z \leq \wtZ$, the inequality \eqref{r14} in this region has been verified in \cite[Section 4.2]{FeKrVa}.
Consequently, it is enough to check \eqref{r14} for $Z \geq \wtZ$. Introducing a new quantity 
$y = \frac{\vr}{\tvr}$ we have to show
\begin{equation} \label{r15}
F(y,Z) = 	S(Z) - S(\wtZ) + \frac{5}{2} - \frac{1}{y}	- \frac{3}{2} \frac{P(Z)}{Z^{\frac{5}{3}}} y^{\frac{2}{3}} \wtZ^{\frac{2}{3}} + 
\frac{1}{10} |S(Z) - S(\wtZ) |^2 \leq 0
\end{equation}
for all $Z \geq \wtZ$, $y > 0$.

To begin, it is also easy to check 
\begin{align}
F(y,Z) &\to - \infty \ \mbox{as}\ y \to 0 \ \mbox{for any fixed}\ Z \geq \widetilde{Z}, \br
F(y,Z) &\to - \infty \ \mbox{as}\ y \to \infty \ \mbox{for any fixed}\ Z \geq \widetilde{Z}.
\label{m4}
\end{align}
	
Next, compute 
\begin{equation} 
\frac{\partial F(y,Z)}{\partial y} = 
\frac{1}{y^2} - \frac{P(Z)}{Z^{\frac{5}{3}}} \widetilde{Z}^{\frac{2}{3}} \frac{1}{y^{\frac{1}{3}}}, 
\label{m5}
\end{equation}		
\begin{align} 
\frac{\partial F(y,Z)}{\partial Z} &= S'(Z) + y^{\frac{2}{3}} 
\left( \frac{\widetilde{Z}}{Z} \right)^{\frac{2}{3}} \frac{3}{2} \frac{  \frac{5}{3}P(Z) - P'(Z) Z   }{Z^{2}} + 
\frac{1}{5} \Big( S(Z) - S(\widetilde{Z}) \Big) S'(Z) \br 
&= S'(Z) \left(1     -   y^{\frac{2}{3}} 
\left( \frac{\widetilde{Z}}{Z} \right)^{\frac{2}{3}}                 + \frac{1}{5}   \Big( S(Z) - S(\widetilde{Z}) \Big)              \right), 
\label{m6} 
	\end{align}
and 
\begin{equation} \label{m7}
\frac{\partial^2 F(y,Z)}{\partial y^2} = 
-2 \frac{1}{y^3} + \frac{1}{3} \frac{P(Z)}{Z^{\frac{5}{3}}} \widetilde{Z}^{\frac{2}{3}} \frac{1}{y^{\frac{4}{3}}}, 
\end{equation}
\begin{align} 
	\frac{\partial^2 F(y,Z)}{\partial Z^2} &= S''(Z) \left(1     -   y^{\frac{2}{3}} 
	\left( \frac{\widetilde{Z}}{Z} \right)^{\frac{2}{3}}                 + \frac{1}{5}   \Big( S(Z) - S(\widetilde{Z}) \Big)              \right) \br 
	&+ S'(Z) \left( \frac{1}{5} S'(Z) + \frac{2}{3} y^{\frac{2}{3}}  \frac{\widetilde{Z}^{\frac{2}{3}} }{Z^{\frac{5}{3}}}   \right), \br
\frac{\partial^2 F(y,Z)}{\partial y \partial Z} &= 
\widetilde{Z}^{\frac{2}{3}} \frac{1}{y^{\frac{1}{3}}}	  \frac{  \frac{5}{3} P(Z) - P'(Z) Z }{Z^{\frac{8}{3}}} = - \frac{2}{3}
\left( \frac{\widetilde{Z}}{Z} \right)^{\frac{2}{3}} \frac{1}{y^{\frac{1}{3}}} 	 S'(Z).
	\label{m8} 
\end{align}
Consequently, 
\begin{align} 
\frac{ \partial F(1, \widetilde{Z} )}{\partial y} &= 
\frac{ \partial F(1, \widetilde{Z} )}{\partial Z} = 0, \br
\frac{ \partial^2 F(1, \widetilde{Z}) }{\partial y^2} &= - \frac{5}{3},\ 
\frac{ \partial^2 F(1, \widetilde{Z}) }{\partial Z^2} = 
\frac{1}{5} (S'(\widetilde{Z}))^2 + \frac{2}{3} \frac{1}{\widetilde{Z}} S'(\widetilde{Z}), \ 
\frac{ \partial^2 F(1, \widetilde{Z})}{\partial y \partial Z} = 
- \frac{2}{3} S'(\widetilde{Z}). 	\label{m9}
\end{align}
In addition, we suppose 
\begin{equation} \label{m10}
	S'(\widetilde{Z}) = - \frac{1}{\widetilde{Z}}; 
	\end{equation}
whence
\begin{equation} \label{m11}	 
\frac{ \partial^2 F(1, \widetilde{Z}) }{\partial y^2} = - \frac{5}{3},\ 
\frac{ \partial^2 F(1, \widetilde{Z}) }{\partial Z^2} = 
-\frac{7}{15} \frac{1}{\widetilde{Z}^2},\ 
\frac{ \partial^2 F(1, \widetilde{Z})}{\partial y \partial Z} = 
 \frac{2}{3} \frac{1}{\widetilde{Z}} .
\end{equation}
In particular, the Hessian $\nabla^2 F(1,\widetilde{Z})$ is strictly negatively definite and the function $F$ attains a strict local maximum 
$F(1, \wtZ) = 0$. 

Next, given $Z > \wtZ$ we have 
\begin{equation} \label{m14}
	\frac{\partial F(\Ov{y},Z)}{\partial y} = 
	\frac{1}{\Ov{y}^2} - \frac{P(Z)}{Z^{\frac{5}{3}}} \widetilde{Z}^{\frac{2}{3}} \frac{1}{\Ov{y}^{\frac{1}{3}}} = 0 \ 
	\Rightarrow \ \frac{1}{\Ov{y}^{\frac{5}{3}}} = \frac{P(Z)}{Z^{\frac{5}{3}}} \widetilde{Z}^{\frac{2}{3}}, 
	\end{equation}
therefore, in accordance with \eqref{m7}
\[
\frac{\partial^2 F(\Ov{y},Z)}{\partial y^2} = 
-2 \frac{1}{\Ov{y}^3} + \frac{1}{3} \frac{P(Z)}{Z^{\frac{5}{3}}} \widetilde{Z}^{\frac{2}{3}} \frac{1}{\Ov{y}^{\frac{4}{3}}}
= - 2 \frac{1}{\Ov{y}^3} + \frac{1}{3} \frac{1}{\Ov{y}^3} < 0. 
\]	
Consequently, for any fixed $Z > \wtZ$, the function $y \mapsto F(y,Z)$ attains its strict global maximum at the point 
$\Ov{y}$. More specifically, 
\begin{align} 
	\max_{y > 0} &F(y,Z) = F(\Ov{y}, Z) = S(Z) - S(\widetilde{Z}) + \frac{5}{2} - \frac{1}{\Ov{y}} - \frac{3}{2}
	\frac{P(Z)}{Z^{\frac{5}{3}}}  \Ov{y}^{\frac{2}{3}} \widetilde{Z}^{\frac{2}{3}} + 
	\frac{1}{10}|S(Z) - S(\widetilde{Z})|^2 \br 
	&= S(Z) - S(\widetilde{Z}) + \frac{5}{2} - \frac{5}{2}\left( \frac{P(Z)}{Z^{\frac{5}{3}}} \widetilde{Z}^{\frac{2}{3}} \right)^{\frac{3}{5}}   + 
	\frac{1}{10}|S(Z) - S(\widetilde{Z})|^2. 
	\label{m15}
\end{align}

As EOS for $Z \geq \wtZ$ is given by \eqref{M5}, \eqref{M6}, we deduce
\begin{align} 
\max_{y > 0} F(y,Z) = \frac{\wtZ}{Z} + \frac{3}{2} - \frac{5}{2} \left( \frac{3}{5} + \frac{2}{5} \left( \frac{ \widetilde{Z} }{Z} \right)^{\frac{5}{3}} \right)^{\frac{3}{5}}   + 
\frac{1}{10} \left( \frac{\wtZ}{Z} - 1 \right)^2. 
\label{m20}	
\end{align}
Setting $Y = \frac{\wtZ}{Z}$ we have to evaluate 
\[ 
G(Y) = Y + \frac{3}{2} - \frac{5}{2} \left( \frac{3}{5} + \frac{2}{5} Y^{\frac{5}{3}} \right)^{\frac{3}{5}}   + 
\frac{1}{10} \left( Y - 1 \right)^2, \ 0 \leq Y \leq 1.
\]
\[
G'(Y) = \frac{4}{5} -  \left( \frac{3}{5} + \frac{2}{5} Y^{\frac{5}{3}} \right)^{- \frac{2}{5}}  Y^{\frac{2}{3}}  + \frac{1}{5} Y .
\]
\[
G(0) = \frac{3}{2} + \frac{1}{10} - \frac{5}{2} \left( \frac{3}{5} \right)^{\frac{3}{5}} 
= \frac{5}{2} \left( \frac{3}{5} - \left( \frac{3}{5} \right)^{\frac{3}{5}} \right) + \frac{1}{10} < 0,\ 
G(1) = 0.
\]
Moreover, it follows that 
\[
G''(Y) = \frac{1}{5} + \frac{4}{15} \left( \frac{3}{5} + \frac{2}{5} Y^{\frac{5}{3}} \right)^{- \frac{7}{5}} Y^{\frac{4}{3}} 
-  \frac{2}{3} \left( \frac{3}{5} + \frac{2}{5} Y^{\frac{5}{3}} \right)^{- \frac{2}{5}}  Y^{- \frac{1}{3}} ,  
\]
and hence 
$$
 G''(0) = - \infty,\ G''(1) = \frac{7}{15} - \frac{2}{3} < 0.
$$
In addition, we have 
\[
\left( \frac{3}{5} + \frac{2}{5} Y^{\frac{5}{3}} \right)^{- \frac{7}{5}} Y^{\frac{4}{3}} \leq \left( \frac{3}{5} \right)^{-\frac{7}{5}} \leq 1.08.
\]
Consequently, 
\begin{equation} \label{m21}
G''(Y) \leq \frac{1}{5} + \frac{4}{15} 1.08 - \frac{2}{3} \leq \frac{1}{2} - \frac{2}{3} < 0.
\end{equation}
We conclude that, $G$ is a concave function in $(0,1)$ attaining its maximum at $G(1) = 0$, which yields the desired conclusion \eqref{r14}.

\subsection{Relative energy estimates -- conclusion} 

Going back to \eqref{r3} we obtain 
\begin{align}
	&\left[ \intO{ E_\ep \left(\vre, \vte, \vue \Big| \tvr, \tvt, \tvu \right) } \right]_{t = z}^{t = \tau} \br 
	&+ \int_z^\tau \int_{x_1 = L} E\left( \vre, \vte, \vue \Big| \tvr, \tvt, \tvu \right) \tvu_R \D \sigma_x \dt \br 
	&+ \int_z^\tau \intO{ \frac{\tvt}{\vte} \left( \mathbb{S}_\ep (\vte, \Ds \vue) : \Ds \vue + \ep \frac{\kappa(\vte) |\Grad \vte|^2}{\vte} \right) } \dt \br 
	&\leq  \ep \int_z^\tau \intO{  \left( \frac{\kappa(\vte) \Grad \vte}{\vte} \right) \cdot \Grad \tvt  } \dt +
	\int_z^\tau \intO{ \mathbb{S}_\ep(\vte, \Ds \vue) \Big] : \Ds \tvu } \dt  	+ \omega(\ep,z). 
	\label{r27}
\end{align}
Seeing that the remaining two integrals can be absorbed by the left--hand side modulo an error $\omega(\ep, z)$, we conclude 
\begin{equation} \label{r28}
	\intO{ E_\ep \left(\vre, \vte, \vue \Big| \tvr, \tvt, \tvu \right) (\tau, \cdot) } 
	\leq \intO{ E_\ep \left(\vre, \vte, \vue \Big| \tvr, \tvt, \tvu \right) (z, \cdot) } + \omega(\ep,z) 
\end{equation}
for any $\tau \geq z$ and a.a. $z > 0$.

\section{Proof of Theorem \ref{MT1}}
\label{P}

\begin{proof}
Revisiting formula \eqref{r28}, we have 
\begin{align}
&\intO{ E_\ep \left(\vre, \vte, \vue \Big| \tvr, \tvt, \tvu \right) (z, \cdot) } \br &= 
\intO{ \left( \frac{1}{2} \vre |\vue - \tvu|^2 + \vre e_\ep (\vre, \vte) - \tvt \vre s_\ep (\vre, \vte) \right)(z, \cdot) } \br &-  
\intO{ \left(
\Big( e_\ep(\tvr, \tvt) - \tvt s_\ep(\tvr, \tvt) + \frac{p_\ep(\tvr, \tvt)}{\tvr} \Big) \vre + p_\ep(\tvr, \tvt)\right)(z, \cdot) } \br
&= \intO{ \left( \frac{1}{2} \vre |\vue - \vuB|^2 + \vre e_\ep (\vre, \vte) - \vtB \vre s_\ep (\vre, \vte) \right)(z, \cdot) } \br &-  
\intO{ \left(
\Big( e_\ep(\tvr, \tvt) - \tvt s_\ep(\tvr, \tvt) + \frac{p_\ep(\tvr, \tvt)}{\tvr} \Big) \vre + p_\ep(\tvr, \tvt)\right) (z, \cdot)} \br 
&+ \intO{ \left( \vre \vue \cdot (\vuB - \tvu) + \frac{1}{2} \vre \left( |\tvu|^2 - |\vuB|^2 \right) + 
(\vtB - \tvt) \vre \left( s_\ep (\vre, \vte) - s_\ep(\tvr, \tvt) \right)	\right) (z, \cdot)  } \br 
&+ \intO{ (\vtB - \tvt) \vre s_\ep (\tvr, \tvt) (z, \cdot) }. 
\label{P1}
\end{align}

Now, the crucial observation is that we may choose a superharmonic (convex) function $\vtB$ satisfying \eqref{e3}, and 
\begin{equation} \label{P2} 
0 < \vtB \leq \tvt .
\end{equation}

By virtue of \eqref{e17}, we get
\begin{align}
	&\intO{ E_\ep \left(\vre, \vte, \vue \Big| \tvr, \tvt, \tvu \right) (z, \cdot) } \br 
	&\leq \intO{ \left( \frac{1}{2} \vr_{0,\ep} |\vu_{0,\ep} - \vuB|^2 + \vr_{0,\ep} e_\ep (\vr_{0,\ep}, \vt_{0,\ep}) - \vtB \vr_{0,\ep} s_\ep (\vr_{0,\ep}, \vt_{0,\ep}) \right) } \br &-  
	\intO{ \left(
		\Big( e_\ep(\tvr, \tvt) - \tvt s_\ep(\tvr, \tvt) + \frac{p_\ep(\tvr, \tvt)}{\tvr} \Big) \vre + p_\ep(\tvr, \tvt)\right) (z, \cdot)} \br 
	&+ \intO{ \left( \vre \vue \cdot (\vuB - \tvu) + \frac{1}{2} \vre \left( |\tvu|^2 - |\vuB|^2 \right) + 
		(\vtB - \tvt) \vre \left( s_\ep (\vre, \vte) - s_\ep(\tvr, \tvt) \right)	\right) (z, \cdot)  } \br 
	&+ \intO{ (\vtB - \tvt) \vre s_\ep (\tvr, \tvt) (z, \cdot) }. 
	\label{P3}
\end{align} 
Moreover, in view of the uniform bounds established in Lemma \ref{eL1}, the functions $\vre$, $\vue$ are uniformly weakly continuous in time. Consequently, 
\begin{align}
	&\intO{ E_\ep \left(\vre, \vte, \vue \Big| \tvr, \tvt, \tvu \right) (z, \cdot) } \br 
	&\leq \intO{ \left( \frac{1}{2} \vr_{0,\ep} |\vu_{0,\ep} - \vuB|^2 + \vr_{0,\ep} e_\ep (\vr_{0,\ep}, \vt_{0,\ep}) - \vtB \vr_{0,\ep} s_\ep (\vr_{0,\ep}, \vt_{0,\ep}) \right) } \br &-  
	\intO{ \left(
		\Big( e_\ep(\tvr, \tvt) - \tvt s_\ep(\tvr, \tvt) + \frac{p_\ep(\tvr, \tvt)}{\tvr} \Big)(0, \cdot) \vr_{0,\ep} + p_\ep(\tvr, \tvt) (0, \cdot) \right) } \br 
	&+ \intO{ \left( \vr_{0,\ep} \vu_{0,\ep} \cdot (\vuB - \tvu(0, \cdot)) + \frac{1}{2} \vr_{0,\ep} \left( |\tvu (0, \cdot) |^2 - |\vuB|^2 \right) \right)} \br &+ 
	\intO{	(\vtB - \tvt) \vre \left( s_\ep (\vre, \vte) - s_\ep(\tvr, \tvt) \right) (z, \cdot)  } \br 
	&+ \intO{ (\vtB - \tvt) \vr_{0,\ep} s_\ep (\tvr, \tvt) (0, \cdot) } + \delta(z), 
	\label{P4}
\end{align} 
where 
\[
\delta(z) \to 0 \ \mbox{as}\ z \to 0 
\]
uniformly in $\ep$.
Finally, it follows from the entropy inequality \eqref{WS4},
\begin{equation} \label{P5}
	\liminf_{z \to 0} \intO{ (\tvt - \vtB) \vre s_\ep (\vre, \vte) } \geq \intO{ (\tvt(0, \cdot) - \vtB) \vr_{0,\ep} s_{\ep} (\vr_{0, \ep}, \vt_{0, \ep}) }
\end{equation}
uniformly in $\ep$, whence we may rewrite \eqref{P4} as 
\begin{align}
	&\intO{ E_\ep \left(\vre, \vte, \vue \Big| \tvr, \tvt, \tvu \right) (z, \cdot) } \br 
	&\leq \intO{ \left( \frac{1}{2} \vr_{0,\ep} |\vu_{0,\ep} - \tvu(0, \cdot)|^2 + \vr_{0,\ep} e_\ep (\vr_{0,\ep}, \vt_{0,\ep}) - \tvt(0, \cdot) \vr_{0,\ep} s_\ep (\vr_{0,\ep}, \vt_{0,\ep}) \right) } \br &-  
	\intO{ \left(
		\Big( e_\ep(\tvr, \tvt) - \tvt s_\ep(\tvr, \tvt) + \frac{p_\ep(\tvr, \tvt)}{\tvr} \Big)(0, \cdot) \vr_{0,\ep} + p_\ep(\tvr, \tvt) (0, \cdot) \right) }  + \delta(z) \br
 &= \intO{ E_\ep \left(\vr_{0,\ep}, \vt_{0,\ep}, \vu_{0,\ep} \Big| \tvr(0, \cdot) , \tvt (0,\cdot), \tvu(0, \cdot  \right) } + \delta(z). 
	\nonumber
\end{align} 

Thus going back to \eqref{r28} we may infer that 
\[
	\intO{ E_\ep \left(\vre, \vte, \vue \Big| \tvr, \tvt, \tvu \right) (\tau, \cdot) } 
\leq \intO{ E_\ep \left(\vr_{0,\ep}, \vt_{0,\ep}, \vu_{0,\ep} \Big| \tvr(0, \cdot) , \tvt (0,\cdot), \tvu(0, \cdot  \right) } + \delta(z) + \omega(\ep,z).
\]
Therefore, chosing first $z> 0$ small enough and then letting $\ep \to 0$ yields the desired conclusion. 
We have proved Theorem \ref{MT1}. 
\end{proof}

\section*{Data availability statement}
Data sharing is not applicable to this article as no data sets were generated or analysed during the current study.

 \section*{Conflict of Interest}
There is no conflict of interests.

\def\cprime{$'$} \def\ocirc#1{\ifmmode\setbox0=\hbox{$#1$}\dimen0=\ht0
	\advance\dimen0 by1pt\rlap{\hbox to\wd0{\hss\raise\dimen0
			\hbox{\hskip.2em$\scriptscriptstyle\circ$}\hss}}#1\else {\accent"17 #1}\fi}


\begin{thebibliography}{10}
	
	\bibitem{CHHS}
	T.~Chang and L.~Hsiao.
	\newblock {\em The {R}iemann problem and interaction of waves in gas dynamics}.
	\newblock Longman Sci. and Tech., Essex, 1989.
	
	\bibitem{ChauFei}
	N.~Chaudhuri and E.~Feireisl.
	\newblock Navier-{S}tokes-{F}ourier system with {D}irichlet boundary
	conditions.
	\newblock {\em Appl. Anal.}, {\bf 101}(12):4076--4094, 2022.
	
	\bibitem{CheChe}
	G.-Q. Chen and J.~Chen.
	\newblock Stability of rarefaction waves and vacuum states for the
	multidimensional {E}uler equations.
	\newblock {\em J. Hyperbolic Differ. Equ.}, {\bf 4}(1):105--122, 2007.
	
	\bibitem{ChenPer1}
	G.-Q. Chen and M.~Perepelitsa.
	\newblock Vanishing viscosity limit of the {N}avier-{S}tokes equations to the
	{E}uler equations for compressible fluid flow.
	\newblock {\em Comm. Pure Appl. Math.}, {\bf 63}(11):1469--1504, 2010.
	
	\bibitem{FeKr2014}
	E.~Feireisl and O.~Kreml.
	\newblock Uniqueness of rarefaction waves in multidimensional compressible
	{E}uler system.
	\newblock {\em J. Hyperbolic Differ. Equ.}, {\bf 12}(3):489--499, 2015.
	
	\bibitem{FeKrVa}
	E.~Feireisl, O.~Kreml, and A.~Vasseur.
	\newblock Stability of the isentropic {R}iemann solutions of the full
	multidimensional {E}uler system.
	\newblock {\em SIAM J. Math. Anal.}, {\bf 47}(3):2416--2425, 2015.
	
	\bibitem{FeNo6A}
	E.~Feireisl and A.~Novotn\'y.
	\newblock {\em Singular limits in thermodynamics of viscous fluids}.
	\newblock Advances in Mathematical Fluid Mechanics. Birkh\"auser/Springer,
	Cham, 2017.
	\newblock Second edition.
	
	\bibitem{FeiNovOpen}
	E.~Feireisl and A.~Novotn{\' y}.
	\newblock {\em Mathematics of open fluid systems}.
	\newblock Birkh{\" a}user--Verlag, Basel, 2022.
	
	\bibitem{GooXin}
	J.~Goodman and Z.~P. Xin.
	\newblock Viscous limits for piecewise smooth solutions to systems of
	conservation laws.
	\newblock {\em Arch. Rational Mech. Anal.}, {\bf 121}(3):235--265, 1992.
	
	\bibitem{HofLiu}
	D.~Hoff and T.-P. Liu.
	\newblock The inviscid limit for the {N}avier-{S}tokes equations of
	compressible, isentropic flow with shock data.
	\newblock {\em Indiana Univ. Math. J.}, {\bf 38}(4):861--915, 1989.
	
	\bibitem{LiHWanWan}
	H.-L. Li, T.~Wang, and Y.~Wang.
	\newblock Wave phenomena to the three-dimensional fluid-particle model.
	\newblock {\em Arch. Ration. Mech. Anal.}, {\bf 243}(2):1019--1089, 2022.
	
	\bibitem{LiWaWa}
	L.-A. Li, D.~Wang, and Y.~Wang.
	\newblock Vanishing viscosity limit to the planar rarefaction wave for the
	two-dimensional compressible {N}avier-{S}tokes equations.
	\newblock {\em Comm. Math. Phys.}, {\bf 376}(1):353--384, 2020.
	
	\bibitem{LiWanWan}
	L.-A. Li, D.~Wang, and Y.~Wang.
	\newblock Vanishing dissipation limit to the planar rarefaction wave for the
	three-dimensional compressible {N}avier-{S}tokes-{F}ourier equations.
	\newblock {\em J. Funct. Anal.}, {\bf 283}(2):Paper No. 109499, 50, 2022.
	
	\bibitem{LiWanTWan}
	L.-A. Li, T.~Wang, and Y.~Wang.
	\newblock Stability of planar rarefaction wave to 3{D} full compressible
	{N}avier-{S}tokes equations.
	\newblock {\em Arch. Ration. Mech. Anal.}, {\bf 230}(3):911--937, 2018.
	
	\bibitem{LiLi1}
	X.~Li and L.-A. Li.
	\newblock Vanishing viscosity limit to the planar rarefaction wave for the
	two-dimensional full compressible {N}avier-{S}tokes equations.
	\newblock {\em J. Differential Equations}, {\bf 269}(4):3160--3195, 2020.
	
	\bibitem{Xin93}
	Z.~P. Xin.
	\newblock Zero dissipation limit to rarefaction waves for the one-dimensional
	{N}avier-{S}tokes equations of compressible isentropic gases.
	\newblock {\em Comm. Pure Appl. Math.}, 46(5):621--665, {\bf 1993}.
	
\end{thebibliography}


\end{document}